\documentclass{amsart}
\usepackage[english]{babel}
\usepackage[latin1]{inputenc}
\usepackage[dvips,final]{graphics}
\usepackage{amsmath,amsfonts,amssymb,amsthm,amscd,array,
stmaryrd,mathrsfs,mathdots,epigraph}
\usepackage[makeroom]{cancel}
\usepackage{pstricks}
\usepackage[all]{xy}
\usepackage{url}
\usepackage{multirow, blkarray}
\usepackage{booktabs}
\usepackage{caption}
\usepackage{subcaption}

\usepackage{blindtext}
\usepackage{textcomp}
\usepackage[final]{epsfig}
\usepackage{color}
\usepackage{mathabx}

\theoremstyle{plain}
\newtheorem{pb}{Problem}%[section]
\newtheorem{thm}{Theorem}%[section]

\theoremstyle{definition}
\newtheorem*{defn}{Definition}

\newtheorem{exe}{Exercise}
\newtheorem*{ex}{Example}
\newtheorem*{hint}{Hint}

\newcommand{\Z}{\mathbb{Z}}

\newcommand{\Id}{\mathrm{Id}}
\newcommand{\SL}{\mathrm{SL}}
\newcommand{\PSL}{\mathrm{PSL}}

\def\thup{\mathop{\rm th}\nolimits}

\begin{document}

\title{Counting quiddities of polygon dissections}

\author{Charles H.\ Conley}
\address{
Charles H.\ Conley,
Department of Mathematics 
\\University of North Texas 
\\Denton TX 76203, USA} 
\email{conley@unt.edu}

\author{Valentin Ovsienko}
\address{
Valentin Ovsienko,
Centre National de la Recherche Scientifique,
Laboratoire de Math\'ematiques de Reims, UMR9008 CNRS,
Universit\'e de Reims Champagne-Ardenne,
U.F.R. Sciences Exactes et Naturelles,
Moulin de la Housse - BP 1039,
51687 Reims cedex 2,
France}
\email{valentin.ovsienko@univ-reims.fr}

\thanks{\noindent
C.H.C.\ was partially supported by Simons Foundation Collaboration Grant~519533.\\
\indent V.O.\ was partially supported by the ANR project PhyMath, ANR-19-CE40-0021.
}

\maketitle

\thispagestyle{empty}

A \textit{dissection} is a partition of a convex polygon 
into sub-polygons, or \textit{cells}, by non-crossing diagonals.
The simplest and most popular class
of dissections is the \textit{triangulations},
where the cells are all triangles.
Polygon dissections have appealed to mathematicians
since the time of ancient Greece,
and they play roles in the study of such
modern objects as associahedra and cluster algebras.
The entry A033282 in the On-Line Encyclopedia of Integer Sequences
(OEIS)~\cite{OEIS} gives an extensive list of references.

The enumeration of dissections was initiated by Euler~\cite{Eul},
who already knew of the Catalan numbers and their main properties.
The topic arises in several areas of mathematics,
connecting Young tableaux, knot invariants,
continued fractions, and many other notions
(some discussed at the end of this article).
It is still an active area; in particular, the enumeration of classes
of dissections under various equivalence relations
has been considered by a number of modern authors.
In this brief note we advertise a new open problem of this nature.

%%%%%%%%%%%%%%%%%%%%%%%%
\section*{Classical formulas}\label{ClaS}
%%%%%%%%%%%%%%%%%%%%%%%%

Here are the most classical sequences arising in our context:

\begin{itemize}

\item
The \textit{Catalan sequence} [OEIS~A000108] is
\begin{equation} 
\label{C_n}
   C_n := \frac{1}{n+1}\binom{2n}{n},
\end{equation}
the number of triangulations of the $(n+2)$-gon.
For hundreds of interpretations and references,
see Stanley's addendum~\cite{StaC} and Pak's webpage \cite{Pak}.

\smallbreak \item
The \textit{Kirkman-Cayley sequence} [OEIS~A033282] is
\begin{equation} 
\label{KC}
   D_{n,m} := \frac{1}{n+1} \binom{n-1}{m-1} \binom{n+m}{m},
\end{equation}
the number of dissections of the $(n+2)$-gon into $m$ cells.
It contains the Catalan sequence as $D_{n,n}$.

\end{itemize}

The Kirkman-Cayley formula~\eqref{KC} was
conjectured by Kirkman in 1857 \cite{Kir}
and proved by Cayley in 1891 \cite{Cay}.
It was also stated as a question by Prouhet in 1866 \cite{Pro}.
Many modern proofs are available;
to give only a few examples,
short proofs using generating functions
may be found in \cite{Rea} and~\cite{FN},
and proofs by bijection are given in
\cite{Sta} (using Young tableaux),
\cite{PS} (in relation to knot theory), and~\cite{Gai}.
Later in this article we will explain the generating function approach.

Let us mention a generalization of the Catalan sequence,
the \textit{Fuss sequence}
\begin{equation*}
   F_{n, m} := \frac{1}{n+1} \binom{n+m}{m},
\end{equation*}
the number of dissections of the $(n+2)$-gon into $m$ cells of equal size.
Here $n$ must be a multiple of $m$, and the cells are
all $(\frac{n}{m}+2)$-gons.  In light of Exercise~\ref{Fuss},
the Fuss sequence is not directly related to our subject.

%%%%%%%%%%%%%%%%%%%%%%%%
\section*{Conway-Coxeter quiddities} \label{CoCoS}
%%%%%%%%%%%%%%%%%%%%%%%%

The most obvious equivalence classes of dissections arise
from the action of the dihedral group by rotations and reflections.
Their enumeration has received considerable attention;
see for example \cite{BR, Rea}.

Here we propose the enumeration of the classes of a
different type of equivalence, given by the notion of \textit{quiddity}.
The idea goes back to Conway and Coxeter,
who used quiddities of triangulations to classify
\textit{frieze patterns} \cite{CoCo}
(see also \cite{Bau} and the video~\cite{Tab}).
Their definition applies equally to arbitrary dissections:

\begin{defn}
The \textit{quiddity} of a dissection of the $N$-gon
is the $N$-tuple $(c_1, \ldots, c_N)$, 
where $c_i$ is the number of cells contacting
the $i^{\thup}$ vertex.
\end{defn}

Here are some simple examples of dissections, depicted with their quiddities:
$$
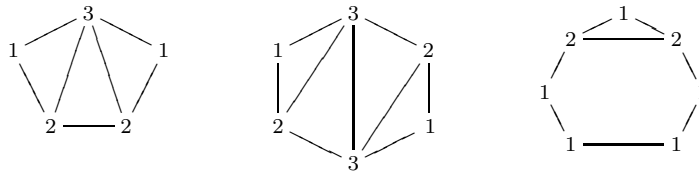
\begin{footnotesize}
\xymatrix @!0 @R=0.50cm @C=0.5cm
{
&& 3 \ar@{-}[rrd]\ar@{-}[lld]\ar@{-}[lddd]\ar@{-}[rddd]&
\\
1 \ar@{-}[rdd]&&&& 1 \ar@{-}[ldd]\\
\\
& 2 \ar@{-}[rr]&& 2
}
\qquad \qquad
\xymatrix @!0 @R=0.50cm @C=0.5cm
{
&&3\ar@{-}[rrd]\ar@{-}[lld]\ar@{-}[llddd]\ar@{-}[dddd]&
\\
1\ar@{-}[dd]&&&& 2\ar@{-}[llddd]\ar@{-}[dd]\\
\\
2\ar@{-}[rrd]&&&& \,1\ar@{-}[lld]\\
&&3
}
\qquad \qquad
\xymatrix @!0 @R=0.35cm @C=0.35cm
{
&&&1\ar@{-}[rrd]\ar@{-}[lld]
\\
&2\ar@{-}[ldd]\ar@{-}[rrrr]&&&& 2\ar@{-}[rdd]&\\
\\
1\ar@{-}[rdd]&&&&&& 1\ar@{-}[ldd]\\
\\
&1\ar@{-}[rrrr]&&&&1
}
\end{footnotesize}
$$
The quiddity gives a good deal of information about the dissection.
For example, consider the following exercises:

\begin{exe}
The quiddity sum $c_1 + \cdots + c_N$
determines the total number of cells in the dissection.
\end{exe}

\begin{hint}
Relate both the sum and the number of cells to
the total number of diagonals of the polygon in the dissection.
\end{hint}

\begin{exe} \label{Fuss}
Dissections into cells of equal sizes are determined by their quiddities.
\end{exe}

\begin{hint}
In any dissection, at least one cell is ``exterior'':
only one of its sides is a diagonal of the polygon.
Check that the quiddity determines the locations of the external cells,
and then base an inductive argument
on the operation of removing an exterior cell.
\end{hint}

\begin{figure}[h!]
\begin{footnotesize}
$$
\xymatrix @!0 @R=0.30cm @C=0.5cm
 {
&1\ar@{-}[ldd]\ar@{-}[rr]&& 2\ar@{-}[rdd]\\
\\
2\ar@{-}[dd]\ar@{-}[rrruu]&&&&1\ar@{-}[dd]\\
\\
1\ar@{-}[rdd]&&&&2\ar@{-}[ldd]\ar@{-}[llldd]\\
\\
&2\ar@{-}[rr]&& 1
}
\qquad\qquad\qquad
\xymatrix @!0 @R=0.30cm @C=0.5cm
 {
&1\ar@{-}[ldd]\ar@{-}[rr]&& 2\ar@{-}[rdd]\\
\\
2\ar@{-}[dd]\ar@{-}[rdddd]&&&&1\ar@{-}[dd]\\
\\
1\ar@{-}[rdd]&&&&2\ar@{-}[ldd]\ar@{-}[luuuu]\\
\\
&2\ar@{-}[rr]&& 1
}
$$
\caption{Distinct dissections with the same quiddity}
\label{ExFig}
\end{footnotesize}
\end{figure}
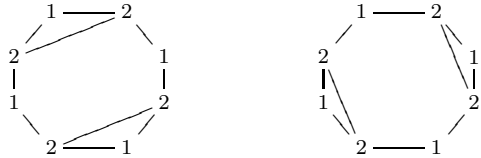

Exercise~\ref{Fuss} shows in particular that
triangulations are determined by their quiddities.
This is not true for arbitrary dissections.
As can be seen in Figure~\ref{ExFig},
it begins to fail in the octagonal case.
The two dissections in the figure are congruent by rotation,
but this is not always the case:

\begin{exe}
Construct distinct dissections with the same quiddity
which are not congruent by any dihedral symmetry.
\end{exe}

We formulate the following general problem.
To the best of our knowledge, it is open and
has not previously been considered.
It seems (at least to us!) to be difficult; at any rate,
more difficult than enumerating the dissections themselves.

\begin{pb}
\label{ThePb}
Enumerate the distinct quiddities of dissections
of the $N$-gon into $m$ cells.
\end{pb}

%%%%%%%%%%%%%%%%%%%%%%%%
\section*{$3$-periodic quiddities}
%%%%%%%%%%%%%%%%%%%%%%%%

Here we give the enumeration of the quiddities of a particular class of dissections.

\begin{defn}
An \textit{$\ell$-periodic dissection} of a polygon
is a dissection such that the number of vertices
of every cell is $3$ modulo~$\ell$.
\end{defn}

For example, $1$-periodic dissections are simply arbitrary dissections,
$2$-periodic dissections have only cells with odd numbers of sides,
and in $3$-periodic dissections,
the number of vertices of every cell is a multiple of~$3$.
The dissections depicted in the preceding section are all $3$-periodic.
Our result is as follows.

\begin{thm}[\cite{CVArx}] \label{Ourtheorem}
Let $Q_{n, m}$ be the number of distinct quiddities
of $3$-periodic dissections of the $(n+2)$-gon into $m$ cells.
Then $Q_{n, m} = 0$ unless $n \equiv m$ mod~$3$,
in which case it is
\begin{equation} \label{Ourformula}
   Q_{n, m} = \sum_{0 \le s \le (n-m)/3}
   \frac{n-m-3s+2}{n-s+1}\, \binom{m+s -2}{s} \binom{n+m-s-1}{m-1}.
\end{equation}
\end{thm}

The proof is somewhat involved
and relies heavily on $3$-periodicity;
we give a little of its flavor in the next two sections.
As far as we can see, it does not adapt to
the case of arbitrary dissections.
Here are the initial values of $Q_{n, m}$:

\begin{footnotesize}
\begin{center}
\begin{tabular}{ |c||c|c|c|c|c|c|c|c|c|c|c|c|c|c|c|}
\hline
$\;\;{}_{\textstyle m} \backslash{} {\textstyle n}$ 
&0 & 1 & 2 & 3 & 4 & 5 & 6 & 7 & 8 & 9 & 10 &11 & 12 &13&14\\
\hline\hline
$n$ & 1 & 1 & 2 & 5 & 14 & 42 & 132 & 429 & 1430 & 4862 & 16796 & 58786 & 208012 & 742900 & 2674440 \\ 
\hline
$n-3$ &  &&&& 1 & 7 & 34 & 147 & 605 & 2431 & 9646 & 38012 & 149226 & 584630 & 2288132\\ 
\hline
$n-6$ &  & &&&&&& 1& 15 & 121 & 758 & 4160 & 21098 & 101660 & 472872\\ 
\hline
$n-9$ &&&&&&&&&&& 1& 26 & 315 & 2710 & 19234\\ 
\hline
$n-12$ &&&&&&&&&&&&&& 1 & 40\\ 
\hline
\end{tabular}
\medbreak
\captionof{table}{The coefficients $Q_{n, m}$}
\end{center}
\end{footnotesize}
Note that the first row, $Q_{n, n}$, is the Catalan sequence, and
the second row, $Q_{n, n-3}$, counts quiddities of dissections
of the $(n+2)$-gon into $n-4$ triangles and $1$ hexagon.
The column sums are the total number of quiddities
of $3$-periodic dissections of the $(n+2)$-gon
[OEIS A348666].  

In the last section we will see that $3$-periodic dissections
have applications to other areas of mathematics.
Surprisingly, they are also the only class of dissections
whose quiddities we have been able to count.

%%%%%%%%%%%%%%%%%%%%%%%%
\section*{Generating functions}\label{GFS}
%%%%%%%%%%%%%%%%%%%%%%%%

It is always useful to organize sequences of integers as power series.
Such series are called \textit{generating functions},
and their analytic properties are powerful tools.
For example, many classes of dissections may be enumerated as follows:
find a functional equation for the generating function,
and then apply \textit{Lagrange inversion}.
In this section we outline the derivation of \eqref{C_n},
\eqref{KC}, and \eqref{Ourformula} via this strategy.

%%%%%%%%
\subsection*{The case of Catalan}
%%%%%%%%

Here the generating function is
\begin{equation*}
   C(z) := \sum_{n = 0}^\infty C_n z^n,
\end{equation*}
the coefficient $C_n$ being the number of triangulations of the $(n+2)$-gon.
Wikipedia gives six proofs of the formula \eqref{C_n} for $C_n$,
including some using paths on grids, some using Dyck words,
and one using nothing but ingenious markings of triangulations.

The most basic approach rests on the recurrence relation
\begin{equation*}
   C_n = C_0 C_{n-1} + C_1 C_{n-2} + 
   C_2 C_{n-3} + \cdots + C_{n-1} C_0,
\end{equation*}
where $C_0$ is defined to be~$1$.
To prove it, designate arbitrarily
one edge of the $(n+2)$-gon to be its \textit{base edge},
a standard trick in this field.
There are $n$ triangular cells containing this base edge.
Choose one of them, as in the figure.
It splits the polygon into two pieces,
an $(n_0 + 2)$-gon and an $(n_1 + 2)$-gon,
where $n_0 + n_1 = n-1$.
Therefore there are $C_{n_0} C_{n_1}$
triangulations containing it.
Summing over the $n$ choices gives the relation.
$$
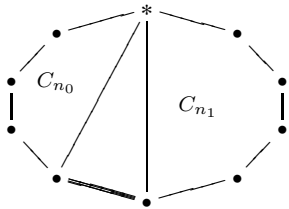
\begin{footnotesize}
\xymatrix @!0 @R=0.32cm @C=0.6cm
 {
&&&\ast\ar@{-}[dddddddd]\ar@{-}[llddddddd]\ar@{-}[lld]\ar@{-}[rrd]&
\\
&\bullet\ar@{-}[ldd]&&&& \bullet\ar@{-}[rdd]\\
\\
\bullet\ar@{-}[dd]&C_{n_0}&&&&&\bullet\ar@{-}[dd]\\
&&&&\;\;C_{n_1}
\\
\bullet\ar@{-}[rdd]&&&&&&\bullet\ar@{-}[ldd]\\
\\
&\bullet\ar@<0.1mm>@{-}@[blue][rrd]\ar@<-0.1mm>@{-}@[blue][rrd]\ar@<-0.2mm>@{-}@[blue][rrd]\ar@<0.2mm>@{-}@[blue][rrd]
&&&& \bullet\ar@{-}[lld]\\
&&&\bullet&
}
\end{footnotesize}
$$

The recurrence relation is equivalent to the functional equation
\begin{equation*}
   C(z) = 1 + z C(z)^2
\end{equation*}
for the generating function. To understand this,
it suffices to observe that for $n > 0$, the right hand side
of the recurrence relation is the coefficient of $z^{n-1}$ in $C(z)^2$.

The functional equation is a quadratic in $C(z)$, whose solution is
\begin{equation*}
   C(z) = \frac{1 - \sqrt{1 - 4z}}{2z}.
\end{equation*}
Computing the Taylor expansion of this formula gives~\eqref{C_n}.

%%%%%%%%
\subsection*{Kirkman-Cayley}
%%%%%%%%

In this case the generating function is \textit{bivariate}: it is
\begin{equation*}
   D(z, w) := \sum_{n, m = 0}^\infty D_{n,m} z^n w^m,
\end{equation*}
where $D_{n,m}$ is the number of
dissections of the $(n+2)$-gon into $m$ cells.
It satisfies the functional equation
\begin{equation} \label{KCFE}
   D = 1 + \frac{w z D^2}{1 - z D} =
   1 + w z D^2 + w z^2 D^3 + w z^3 D^4 + \cdots.
\end{equation}

To understand why, designate again a base edge,
and choose an $(r+2)$-cell containing it.
$$
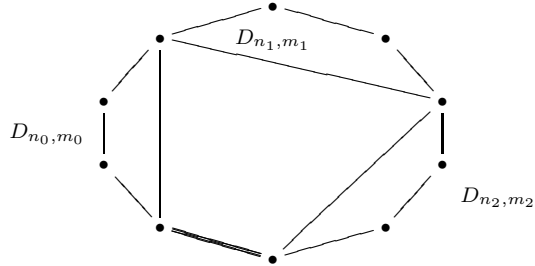
\begin{footnotesize}
\xymatrix @!0 @R=0.42cm @C=0.75cm
{
&&&& \bullet \ar@{-}[lld] \ar@{-}[rrd] &
\\
&& \bullet \ar@{-}[dddddd] \ar@{-}[rrrrrdd] \ar@{-}[ldd]
&& D_{n_1, m_1} && \bullet \ar@{-}[rdd]
\\ \\
& \bullet \ar@{-}[dd] &&&&&& \bullet \ar@{-}[dd] \ar@{-}[lllddddd]
\\
D_{n_0, m_0}
\\
& \bullet \ar@{-}[rdd] &&&&&& \bullet \ar@{-}[ldd]
\\
&&&&&&&& D_{n_2,m_2}
\\
&& \bullet \ar@<0.1mm>@{-}@[blue][rrd] \ar@<-0.1mm>@{-}@[blue][rrd]
\ar@<-0.2mm>@{-}@[blue][rrd] \ar@<0.2mm>@{-}@[blue][rrd]
&&&& \bullet \ar@{-}[lld]
\\
&&&& \bullet &
}
\end{footnotesize}
$$

First question: how many dissections into $m$ cells
contain the chosen $(r+2)$-cell?
The cell splits the polygon into $(r+1)$ pieces:
an $(n_0 + 2)$-gon, an $(n_1 + 2)$-gon, and so on,
up to an $(n_r + 2)$-gon, where $\sum_{s=0}^r n_s = n - r$
(the case $r=2$ is depicted).
Each of these pieces may itself be dissected, and
if the $(n_s + 2)$-gon is dissected into $m_s$ cells,
then we must have $\sum_{s=0}^r m_s = m - 1$.
Therefore the answer is
\begin{equation*}
   \sum_{m_0 + \cdots + m_r = m - 1}
   D_{n_0, m_0} \cdots D_{n_r, m_r}.
\end{equation*}
(Care must be taken with the degenerate case in which
the chosen cell shares some edges with the polygon.
Here $n_s = 0$ for some values of $s$, and
it is necessary to use the convention that
$D_{0, 0} = 1$ and $D_{0, m} = 0$ for $m > 0$.
To put it in words, ``the $2$-gon has a unique dissection,
which has $0$ cells''.)

Second question: in how many dissections into $m$ cells
is the base edge in \textit{some} $(r+2)$-cell?
The set of all $(r+2)$-cells containing the base is indexed
by the $(r+1)$-tuples $(n_0, \ldots, n_r)$ such that
$\sum_{s=0}^r n_s = n - r$, so the answer is
\begin{equation}
\label{rnm}
   \sum_{\substack{n_0 + \cdots + n_r = n - r \\
   m_0 + \cdots + m_r = m - 1}}
   D_{n_0, m_0} \cdots D_{n_r, m_r}.
\end{equation}

Third question: what is the total number $D_{n, m}$
of dissections into $m$ cells?
The base edge may be in a cell of any size, so the answer
is the sum of \eqref{rnm} over~$r$:
\begin{equation*}
   D_{n, m} = \sum_{r = 1}^n \
   \sum_{\substack{n_0 + \cdots + n_r = n - r \\
   m_0 + \cdots + m_r = m - 1}}
   D_{n_0, m_0} \cdots D_{n_r, m_r}.
\end{equation*}

Last question: how does this yield the recurrence relation \eqref{KCFE}?
It suffices to verify that \eqref{rnm} is the coefficient of
$z^n w^m$ in $w z^r D^{r+1}$, or in other words,
the coefficient of $z^{n-r} w^{m-1}$ in $D^{r+1}$.
We leave this as an exercise for the reader.

%%%%%%%%
\subsection*{Lagrange inversion}
%%%%%%%%

The strategy used in the Catalan case to obtain the explicit formula
for $C_n$ from the functional equation for $C(z)$ does not work well
for dissections: one can easily solve for $D$ in \eqref{KCFE},
but it is difficult to use Taylor expansion to obtain~\eqref{KC}.
However, Lagrange's inversion theorem makes things easy.  
We briefly recall it: suppose $\phi(y)$ is a series in $y$
with non-zero constant term.  The theorem gives the terms
of the series $y(z)$ in $z$ which inverts the function
$y \mapsto y / \phi(y)$, i.e., which satisfies $z = y(z) / \phi(y(z))$.
Using the standard notation $[z^n] y$ for the coefficient
of $z^n$ in $y(z)$, the result is
\begin{equation*}
  n [z^n] y = [y^{n-1}] \phi^n.
\end{equation*}

\begin{exe}
Use Lagrange inversion to prove the Kirkman-Cayley formula~\eqref{KC}.
\end{exe}

\begin{hint}
Set $y := zD(z, w)$ and rearrange \eqref{KCFE} to see that one may take
\begin{equation*}
   \phi(y) = 1 + wy / (1 - y - wy).
\end{equation*}
Then use $(1-x)^{-(n+1)} = \sum_{m=0}^\infty \binom{n+m}{m} x^m$.
\end{hint}

Let us remark that the proof of the inversion theorem rests on
nothing but clever use of the residue theorem
and integration by substitution.
Writing $\oint_0$ for integration around zero
and substituting $y/\phi(y)$ for $z$ and $dy$ for $y' dz$
at the appropriate moment, the core of the argument is
\begin{equation*}
   n [z^n] y = [z^{n-1}] y' =
   \oint_0 \frac{y' dz}{2 \pi i z^n} =
   \oint_0 \frac{\phi^n(y) dy}{2 \pi i y^n} =
   [y^{n-1}] \phi^n.
\end{equation*}

Despite its short proof, Lagrange inversion is a powerful tool.
As we have mentioned, it permits the enumeration of
a wide variety of classes of dissections, such as for example
the $\ell$-periodic dissections:

\begin{exe} \label{L}
The number of $\ell$-periodic dissections of the $(n+2)$-gon into $m$ cells is
zero unless $n \equiv m$ mod~$\ell$, in which case it is
\begin{equation*}
   L_{n, m} = \frac{1}{n+1} \binom{m-1 + (n-m)/\ell}{m-1} \binom{n+m}{m}.
\end{equation*}
\end{exe}

\begin{hint}
Adapt the proof of \eqref{KCFE} to deduce that
the corresponding bivariate generating function satisfies
\begin{equation*}
   L(z, w) = 1 + \frac {w z L^2}{1 - z^\ell L^\ell},
\end{equation*}
and then apply Lagrange inversion.
\end{hint}

%%%%%%%%
\subsection*{Quiddities}
%%%%%%%%

Enumerating quiddities appears to be significantly more difficult than
enumerating dissections, and as noted earlier,
we have only been able to accomplish it in the $3$-periodic case.
The difficulty rests in finding functional equations.
Indeed, the $3$-periodic generating function
\begin{equation*}
   Q(z, w) := \sum_{n, m \geq 0} Q_{n, m} z^n w^m
\end{equation*}
does not seem to satisfy a functional equation itself,
but it can be given in terms of an auxiliary function
$P(z, w)$ which is defined by such an equation:
\begin{equation*}
   P(z, w) := 1 + \frac {w z P^2}{1 - z^3 P^2}
   = 1 + w z P^2 + w z^4 P^4 + w z^7 P^6 + \cdots.
\end{equation*}
The main ingredient in the proof of Theorem~\ref{Ourtheorem} is
the formula for $Q$ in terms of $P$:
\begin{equation} \label{QFE}
   Q(z, w) = 1 + \frac {w z P^2}{1 - z^3 P^3}
   = 1 + w z P^2 + w z^4 P^5 + w z^7 P^8 + \cdots.
\end{equation}

From \eqref{QFE}, a few preparatory tricks
and a generalization of Lagrange inversion
known as the Lagrange-B\"urmann formula
give \eqref{Ourformula}.

%%%%%%%%%%%%%%%%%%%%%%%%
\section*{Why is $3$-periodicity special?} \label{Discussion}
%%%%%%%%%%%%%%%%%%%%%%%%

Let us now describe some of the ideas behind the proof of \eqref{QFE}, 
as well as some of the obstacles to counting the quiddities
of other classes of dissections.

%%%%%%%%%%%%%%%%%%%%%%%%
\subsection*{Surgery and maximally open dissections}
%%%%%%%%%%%%%%%%%%%%%%%%

\textit{Surgery} is an operation acting on a single cell
of a dissection to produce a new dissection, as follows:
choose two edges of the cell which are both diagonals of the dissection,
and moreover are separated from one another on both sides
by at least two additional edges of the cell.
Remove the two chosen edges from the dissection
and replace them with the two other
line segments having the same endpoints.
For example, surgery transforms each octagonal dissection
in Figure~\ref{ExFig} into the other.

The crucial point is that surgery does not alter the quiddity:
``surgery equivalence classes refine quiddity equivalence classes''.
Our approach to \eqref{QFE} is to prove that
the two equivalence classes are the same:
any two $3$-periodic dissections with the same quiddity
can be transformed into one another by a sequence of surgeries.

We admit only surgeries preserving $3$-periodicity:
\textit{$3$-periodic surgeries}.  
We define a canonical dissection within each
$3$-periodic surgery equivalence class as follows.
As usual, designate one edge of the polygon to be its base.
This endows each cell with its own naturally defined base edge:
the one closest to the base edge of the polygon.
A surgery on a cell is said to be \textit{opening} if it removes the cell's base edge.
A $3$-periodic dissection is said to be \textit{maximally open}
if it admits no $3$-periodic opening surgeries.

In order to transform a $3$-periodic dissection into a maximally open one,
apply repeated $3$-periodic opening surgeries to all of its cells,
beginning with the cells furthest from the base edge of the polygon.
A simple example is shown.

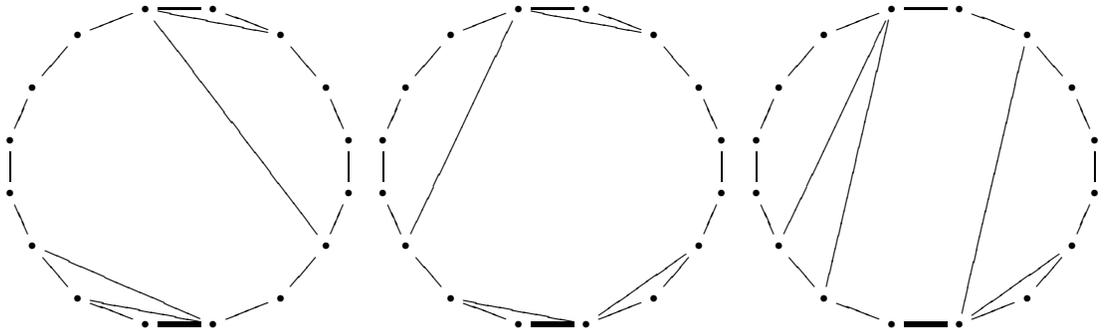
\begin{figure}[htbp]
\begin{tiny}
$$
\xymatrix @!0 @R=0.35cm @C=0.30cm
 {
&&&&&&\bullet\ar@{-}[llld]\ar@{-}[rrr]\ar@{-}[rrrrrrd]
\ar@{-}[rrrrrrrrddddddddd]&&&\bullet\ar@{-}[rrrd]
\\
&&&\bullet\ar@{-}[lldd]&&&&&&&&& \bullet\ar@{-}[rrdd]\\
\\
&\bullet\ar@{-}[ldd]&&&&&&&&&&&&&\bullet\ar@{-}[rdd]\\
\\
\bullet\ar@{-}[dd]&&&&&&&&&&&&&&&\bullet\ar@{-}[dd]\\
&&&&&&&&&&&&\\
\bullet\ar@{-}[rdd]&&&&&&&&&&&&&&&\bullet\ar@{-}[ldd]\\
\\
&\bullet\ar@{-}[rrdd]\ar@{-}[rrrrrrrrddd]&&&&&&&&&&&&&\bullet\ar@{-}[lldd]\\
\\
&&&\bullet\ar@{-}[rrrd]\ar@{-}[rrrrrrd]&&&&&&&&& \bullet\ar@{-}[llld]\\
&&&&&&\bullet\ar@{-}@[blue][rrr]\ar@<0.2mm>@{-}@[blue][rrr]\ar@<-0.2mm>@{-}@[blue][rrr]&&&\bullet
}
\ \
\xymatrix @!0 @R=0.35cm @C=0.30cm
 {
&&&&&&\bullet\ar@{-}[llld]\ar@{-}[rrr]\ar@{-}[rrrrrrd]\ar@{-}[lllllddddddddd]&&&\bullet\ar@{-}[rrrd]
\\
&&&\bullet\ar@{-}[lldd]&&&&&&&&& \bullet\ar@{-}[rrdd]\\
\\
&\bullet\ar@{-}[ldd]&&&&&&&&&&&&&\bullet\ar@{-}[rdd]\\
\\
\bullet\ar@{-}[dd]&&&&&&&&&&&&&&&\bullet\ar@{-}[dd]\\
&&&&&&&&&&&&\\
\bullet\ar@{-}[rdd]&&&&&&&&&&&&&&&\bullet\ar@{-}[ldd]\\
\\
&\bullet\ar@{-}[rrdd]&&&&&&&&&&&&&\bullet\ar@{-}[lldd]\ar@{-}[lllllddd]\\
\\
&&&\bullet\ar@{-}[rrrd]\ar@{-}[rrrrrrd]&&&&&&&&& \bullet\ar@{-}[llld]\\
&&&&&&\bullet\ar@{-}@[blue][rrr]\ar@<0.2mm>@{-}@[blue][rrr]\ar@<-0.2mm>@{-}@[blue][rrr]&&&\bullet
}
\ \
\xymatrix @!0 @R=0.35cm @C=0.30cm
{
&&&&&&\bullet\ar@{-}[llld]\ar@{-}[rrr]\ar@{-}[lllllddddddddd]\ar@{-}[lllddddddddddd]&&&\bullet\ar@{-}[rrrd]
\\
&&&\bullet\ar@{-}[lldd]&&&&&&&&& \bullet\ar@{-}[rrdd]\ar@{-}[lllddddddddddd]\\
\\
&\bullet\ar@{-}[ldd]&&&&&&&&&&&&&\bullet\ar@{-}[rdd]\\
\\
\bullet\ar@{-}[dd]&&&&&&&&&&&&&&&\bullet\ar@{-}[dd]\\
&&&&&&&&&&&&\\
\bullet\ar@{-}[rdd]&&&&&&&&&&&&&&&\bullet\ar@{-}[ldd]\\
\\
&\bullet\ar@{-}[rrdd]&&&&&&&&&&&&&\bullet\ar@{-}[lldd]\ar@{-}[lllllddd]\\
\\
&&&\bullet\ar@{-}[rrrd]&&&&&&&&& \bullet\ar@{-}[llld]\\
&&&&&&\bullet\ar@{-}@[blue][rrr]\ar@<0.2mm>@{-}@[blue][rrr]\ar@<-0.2mm>@{-}@[blue][rrr]&&&\bullet
}
$$
\end{tiny}
\end{figure}

This process leads to the conclusion that there exists
a maximally open dissection in each
$3$-periodic surgery equivalence class.
In order to prove that it is unique,
and simultaneously that $3$-periodic surgery equivalence classes are
the same as quiddity equivalence classes, we prove that
any two maximally open $3$-periodic dissections with the same quiddity
are identical.  This is accomplished by induction on the size of the polygon,
using an argument refining the initial idea of Conway and Coxeter \cite{CoCo}.
The details are given in~\cite{CVArx}.

Thus in the $3$-periodic case, the enumeration of quiddities
is equivalent to the enumeration of maximally open dissections.
From here, we obtain \eqref{QFE} by methods similar to those yielding
the Kirkman-Cayley functional equation~\eqref{KCFE}.

%%%%%%%%%%%%%%%%%%%%%%%%
\subsection*{In search of more general surgery}
%%%%%%%%%%%%%%%%%%%%%%%%

Outside of the $3$-periodic case, surgery equivalence classes
and quiddity equivalence classes are not the same.
For example, consider as a ``toy model'' the
dissections in which all cells are either triangles or quadrilaterals.
Such dissections admit no surgeries, as surgery is never possible
on a cell smaller than a hexagon.
However, we have the following exercise:

\begin{exe}
\label{counterexample}
Dissections into triangles and quadrilaterals
are not determined by their quiddities.
\end{exe}

\begin{hint}
The following dissection has less symmetry than its quiddity:
$$
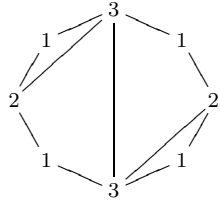
\begin{footnotesize}
\xymatrix @!0 @R=0.40cm @C=0.45cm
{
&&&3\ar@{-}[rrd]\ar@{-}[lld]\ar@{-}[dddddd]\ar@{-}[lllddd]
\\
&1\ar@{-}[ldd]&&&& 1\ar@{-}[rdd]&\\
\\
2\!\ar@{-}[rdd]&&&&&& \!2\ar@{-}[ldd]\ar@{-}[lllddd]\\
\\
&1\ar@{-}[rrd]&&&&1\ar@{-}[lld]\\
&&&3
}
\end{footnotesize}
$$
\end{hint}

Note the similarity to Figure~\ref{ExFig}, which admits a $3$-periodic surgery.
The dissection here can be transformed into its reflection on the vertical
by a more general type of quiddity-preserving
surgery which operates simultaneously on two adjacent cells.
In the next figure one must operate simultaneously on three non-adjacent cells in order to
transform the dissection into its reflection:
$$
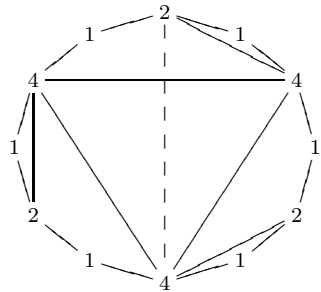
\begin{footnotesize}
\xymatrix @!0 @R=0.30cm @C=0.25cm
{
&&&&&&&&2\ar@{-}[lllld]\ar@{-}[rrrrd]\ar@{-}[rrrrrrrddd]\ar@{--}@[red][ddddddddddd]
\\
&&&&1\ar@{-}[llldd]&&&&&&&& 1\ar@{-}[rrrdd]\\
\\
&4\ar@{-}[rrrrrrrrrrrrrr]\ar@{-}[dddddd]\ar@{-}[lddd]\ar@{-}[rrrrrrrddddddddd]
&&&&&&&&&&&&&&4\ar@{-}[rddd]\ar@{-}[lllllllddddddddd]\\
\\
\\
1\ar@{-}[rddd]&&&&&&&&&&&&&&&&1\ar@{-}[lddd]\\
\\
\\
&2\ar@{-}[rrrdd]&&&&&&&&&&&&&&2\ar@{-}[llldd]\ar@{-}[lllllllddd]\\\
\\
&&&&1&&&&&&&& 1\\
&&&&&&&&4\ar@{-}[llllu]\ar@{-}[rrrru]
}
\end{footnotesize}
$$

Can one define a quiddity-preserving surgery on this type of dissection
whose equivalence classes are equal to the quiddity equivalence classes?
Do such techniques lead to a functional equation?  We do not know.
Let us formulate the relevant simplification of Problem~\ref{ThePb}:

\begin{pb}
\label{PbTwo}
Enumerate the quiddities of dissections into triangles and quadrilaterals.
\end{pb}

Once again, the problem of enumerating the dissections
themselves may be solved using standard techniques:

\begin{exe}
The number of dissections of the $(n+2)$-gon into $m$ cells,
all of which are either triangles or quadrilaterals, is
\begin{equation*}
   D^{3,4}_{n,m} = \frac{1}{n+1} \binom{m}{n-m} \binom{n+m}{m}.
\end{equation*}
\end{exe}

\begin{hint}
As in Exercise~\ref{L}, adapt the proof of \eqref{KCFE} to show
that the bivariate generating function satisfies
$D^{3,4} = 1 + w z (D^{3,4})^2 + w z^2 (D^{3,4})^3$,
and then apply Lagrange inversion.
\end{hint}

Exercise~\ref{counterexample} shows that for arbitrary dissections,
quiddity equivalence classes are larger than
``classical'' surgery equivalence classes.
This is true also in the $2$-periodic case:
here is a dissection into triangles and pentagons with the same quiddity
as its reflection across the vertical:
$$
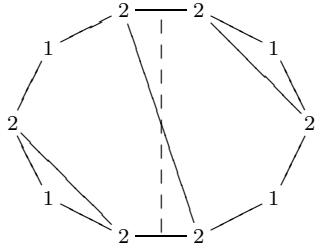
\begin{footnotesize}
\xymatrix @!0 @R=0.50cm @C=0.50cm
{
&&&2\ar@{-}[rrdddddd]\ar@{-}[lld]\ar@{-}[rr]&\ar@{--}@[red][dddddd]&2\ar@{-}[rrd]\ar@{-}[rrrddd]
\\
&1\ar@{-}[ldd]&&&&&& 1\ar@{-}[rdd]&\\
\\
2\!\ar@{-}[rdd]\ar@{-}[rrrddd]&&&&&&&& \!2\ar@{-}[ldd]\\
\\
&1\ar@{-}[rrd]&&&&&&1\ar@{-}[lld]\\
&&&2\ar@{-}[rr]&&2
&&&&&
}
\end{footnotesize}
$$
Is there a type of quiddity-preserving surgery transforming the one into the other?

%%%%%%%%%%%%%%%%%%%%%%%%
\section*{Unexpected connections} \label{Conclusion}
%%%%%%%%%%%%%%%%%%%%%%%%

We conclude with some connections
between quiddities of $3$-periodic dissections and topics in
continued fractions, discrete analysis, and number theory.

%%%%%%%%
\subsection*{Continued fractions}
%%%%%%%%

Every rational number $\frac{r}{s}>1$ has two types
of continued fraction expansions:
$$
\frac{r}{s}
\quad=\quad
a_1 + \cfrac{1}{a_2 
          + \cfrac{1}{\ddots +\cfrac{1}{a_{2m}} } }
           \quad =\quad
c_1 - \cfrac{1}{c_2 
          - \cfrac{1}{\ddots - \cfrac{1}{c_k} } } \,,
$$
where $a_i\geq1$ and $c_i\geq2$.
The second expansion, often called the \textit{Hirzebruch-Jung}
continued fraction, is useful in hyperbolic geometry and toric varieties. 

Conway-Coxeter quiddities appear
in a beautiful manner in this context.
Construct the triangulation in which the integers $a_i$
count the number of consecutive triangles in the same orientation,
alternating base down and base up:
$$
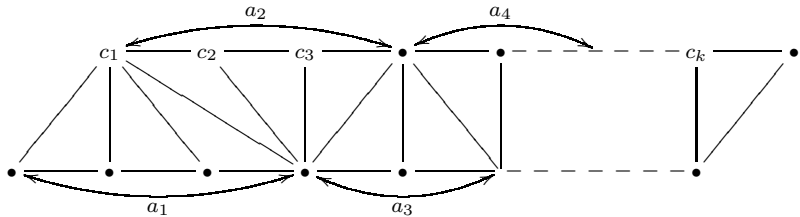
\begin{footnotesize}
\xymatrix @!0 @R=0.8cm @C=1.3cm
{
&c_1\ar@{-}[ldd]\ar@{-}[dd]\ar@{-}[rdd]\ar@{-}[rrdd]
\ar@{-}[r]\ar@/^0.8pc/@{<->}[rrr]^{\textstyle a_2}
&c_2\ar@{-}[r]
\ar@{-}[rdd]&c_3\ar@{-}[r]\ar@{-}[dd]
&\bullet\ar@{-}[r]\ar@{-}[ldd]\ar@{-}[dd]
\ar@{-}[rdd]\ar@{-}[r]\ar@/^0.8pc/@{<->}[rr]^{\textstyle a_4}
&\bullet\ar@{--}[rr]\ar@{-}[dd] &&
c_k\ar@{-}[r]\ar@{-}[dd]& \bullet\ar@{-}[ldd]\\
\\
\bullet\ar@{-}[r]\ar@/_0.8pc/@{<->}[rrr]_{\textstyle a_1}
&\bullet\ar@{-}[r]&\bullet\ar@{-}[r]
&\bullet\ar@{-}[r]\ar@/_0.8pc/@{<->}[rr]_{\textstyle a_3}
&\bullet\ar@{-}[r]&\ar@{--}[rr]&&\bullet&
}
\end{footnotesize}
$$
Remarkably, the integers $c_i$ then turn out to be the
quiddity coefficients on the top of the triangulation!
This should be attributed to~\cite{CoCo};
for an explanation, see~\cite{FB1}.

\begin{ex}
The rational number $\frac{7}{5}$ has the expansions
$$
\frac{7}{5}=
1 + \cfrac{1}{2 
          + \cfrac{1}{1 +\cfrac{1}{1} } }
           \quad =\quad
2 - \cfrac{1}{2 
          - \cfrac{1}{3 } } \,,
$$
corresponding to the triangulation
$$
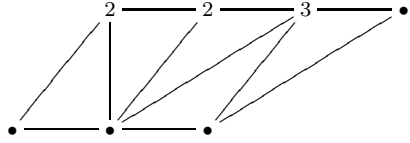
\begin{footnotesize}
\xymatrix @!0 @R=0.8cm @C=1.3cm
{
&2\ar@{-}[ldd]\ar@{-}[dd]\ar@{-}[r]
&2\ar@{-}[ldd]\ar@{-}[r]
&3\ar@{-}[lldd]\ar@{-}[ldd]\ar@{-}[r]
&\bullet\ar@{-}[lldd]\\
\\
\bullet\ar@{-}[r]
&\bullet\ar@{-}[r]
&\bullet
}
\end{footnotesize}
$$
\end{ex}

%%%%%%%%
\subsection*{The discrete Schr\"odinger equation}
%%%%%%%%

Let $(c_i)_{i\in\Z}$ be an $N$-periodic sequence of positive integers:
$c_{i+N}=c_i$ for all $i \in \Z$.
Consider the linear recurrence
\begin{equation}
\label{SL}
v_{i+1}=c_iv_i-v_{i-1},
\end{equation}
where the $v_i$ are unknowns.
This is the simplest second order recurrence,
often called the $1$-dimensional discrete
Schr\"odinger (or Sturm-Liouville) equation.
It appears in many areas of mathematics, and of course
its continuous analog is the subject of a vast literature
(see for example~\cite{Sim}).

\begin{pb}
\label{ThePb2}
Characterize and enumerate the sequences $(c_i)$
such that all solutions $(v_i)$ are either $N$-antiperiodic
($v_{i+N} = -v_i$ for all~$i$) or $N$-periodic.
\end{pb}

%%%%%%%%
\subsection*{The modular group}
%%%%%%%%

Recall that $\SL(2,\Z)$ is the group of
$2\times2$ integer matrices of determinant~$1$, and
$\PSL(2, \Z)$ is its central quotient $\SL(2,\Z)/\{\pm\Id\}$.

\begin{exe}
Every element $A$ of $\SL(2,\Z)$ can be written (not uniquely!)
as a product of ``elementary'' matrices: for some positive integers $c_1, \ldots, c_N$,
\begin{equation}
\label{SLMatEqInt}
A=\begin{pmatrix}
c_1&-1\\[4pt]
1&\phantom{-}0
\end{pmatrix}
\begin{pmatrix}
c_2&-1\\[4pt]
1&\phantom{-}0
\end{pmatrix}
\cdots
\begin{pmatrix}
c_N&-1\\[4pt]
1&\phantom{-}0
\end{pmatrix}.
\end{equation}
\end{exe}

\begin{hint}
Begin with the well-known generators
$
\big(\begin{smallmatrix}
1&1\\[2pt]
0&1
\end{smallmatrix}\big)$
and
$\big(\begin{smallmatrix}
0&-1\\[2pt]
1&\phantom{-}0
\end{smallmatrix}\big)$
of $\SL(2,\Z)$.
\end{hint}

The next problem describes the relations
between elementary matrices in $\PSL(2,\Z)$:

\begin{pb}
\label{ThePb2Bis}
Characterize and enumerate those $c_1, \ldots, c_N$ such that
\eqref{SLMatEqInt} is $\pm\Id$.
\end{pb}

In contrast with Problems~\ref{ThePb} and~\ref{PbTwo}, Problems~\ref{ThePb2}
and~\ref{ThePb2Bis} are in fact already solved:

\begin{thm}[\cite{Val}, Theorem~1.1(i)]
\label{ValRMSthm}
Problems~\ref{ThePb2} and~\ref{ThePb2Bis} have identical solution sets:
the quiddities of the $3$-periodic dissections of the $N$-gon.
\end{thm}

One wonders if the quiddities of other classes
of polygon dissections also have nice interpretations!


\begin{thebibliography}{99}

\bibitem{Bau}
K. Baur, 
{\it Frieze Patterns of Integers,\/}
Math. Intelligencer {\bf 43} (2021), 47--54.

\bibitem{BR}
D.~Bowman, A.~Regev, 
{\it Counting symmetry classes of dissections of a convex regular polygon,\/}
Adv.\ Appl.\ Math.\ {\bf 56} (2014), 35--55.

\bibitem{Cay}
A.~Cayley,
{\it On the partitions of a polygon,\/}
Proc.\ London Math.\ Soc.\ {\bf 22} (1890-1891), 237--262.

\bibitem{CVArx}
C.~Conley, V.~Ovsienko,
{\it Quiddities of polygon dissections and the Conway-Coxeter frieze equation,\/}
arXiv:2107.01234.

\bibitem{CoCo}
J.~H.~Conway, H.~S.~M.~Coxeter,
{\it Triangulated polygons and frieze patterns,\/}
Math.\ Gaz.\ {\bf 57} (1973), 87--94 and 175--183.

\bibitem{Eul}
L.~Euler, 
{\it Letter to Goldbach,\/}
http://eulerarchive.maa.org//correspondence/letters/OO0868.pdf.

\bibitem{FN} 
P.~Flajolet, M.~Noy,
{\it Analytic combinatorics of non-crossing configurations,\/}
Discrete Math.\ {\bf 204} (1999), 203--229.

\bibitem{Gai}
G.~Gaiffi,
{\it Nested sets, set partitions and Kirkman-Cayley dissection numbers,\/}
European J.~Combin.\ {\bf 43} (2016), 279--288.

\bibitem{Kir}
T.~P.~Kirkman,
{\it On the $k$-partitions of the $r$-gon and $r$-ace,\/}
Phil.\ Trans.\ Royal Soc.\ London {\bf 147} (1857), 217--272.

\bibitem{FB1} 
S.~Morier-Genoud, V.~Ovsienko, 
{\it Farey boat:
continued fractions and triangulations,
modular group and polygon dissections,\/}
Jahresber.\ Dtsch.\ Math.-Ver.\ {\bf 121} (2019), no.~2, 91--136.

\bibitem{OEIS} 
The On-Line Encyclopedia of Integer Sequences, OEIS Foundation Inc., http://oeis.org.

\bibitem{Val}
V.~Ovsienko,
{\it Partitions of unity in $\SL(2,\Z)$, negative continued fractions,
and dissections of polygons,\/} 
Res.\ Math.\ Sci.\ {\bf 5} (2018), no.~2, Paper No.~21, 25~pp.

\bibitem{Pak}
I.~Pak,
Catalan numbers page,
https://www.math.ucla.edu/~pak/lectures/Cat/pakcat.htm

\bibitem{Pro}
E.~Prouhet,
{\it Question~774,\/}
Nouvelles Ann.\ Math.\ {\bf 5} (1866), 384.

\bibitem{PS}
J.~Przytycki, A.~Sikora, 
{\it Polygon dissections and Euler, Fuss, Kirkman, and Cayley numbers,\/}
J.~Combin.\ Theory Ser.~A {\bf 92} (2000), 68--76.

\bibitem{Rea}
R.~Read, 
{\it On general dissections of a polygon,\/}
Aequationes Math.\ {\bf 18} (1978), 370--388. 

\bibitem{Sim}
B.~Simon, 
{\it The classical moment problem as a self-adjoint finite difference operator,\/} 
Adv.\ Math.\ {\bf 137} (1998), 82--203. 

\bibitem{Sta}
R.~Stanley,
{\it Polygon dissections and standard Young tableaux,\/}
J.~Combin.\ Theory Ser.~A {\bf 76} (1996), 175--177. 

\bibitem{StaC}
R.~Stanley,
Catalan addendum to Enumerative Combinatorics,
http://www-math.mit.edu/~rstan/ec/catadd.pdf

\bibitem{Tab}
S.~Tabachnikov,
{\it Frieze Patterns} - Numberphile,
https://www.youtube.com/watch?v=0mXz-NP-raY

\end{thebibliography}
\end{document}